\documentclass[10pt]{amsart}
\usepackage{amssymb}
\usepackage{amsmath}
\usepackage[all]{xypic}
\usepackage{epsbox}

\newtheorem{thm}{Theorem}[section]
\newtheorem{lem}[thm]{Lemma}

\newtheorem{prop}[thm]{Proposition}

\newtheorem{kleinthm}{Theorem}[section]

\theoremstyle{definition}

\newtheorem{rem}[thm]{Remark}
\newtheorem{defn}[thm]{Definition}

\theoremstyle{remark}

\numberwithin{equation}{section}

\def\R{{\mathbb R}}
\def\Z{{\mathbb Z}}
\def\C{{\mathbb C}}
\def\Q{{\mathbb Q}}

\def\Gal{\text{\rm Gal}}

\def\Hom{\text{\rm Hom}}
\def\Aut{\text{\rm Aut}}

\def\G{\Gamma}

\allowdisplaybreaks[1]

\begin{document}
\title[A nontrivial algebraic cycle in
the Jacobian variety]
{A nontrivial algebraic cycle in the\\
Jacobian variety of the Klein quartic}

\author[Yuuki Tadokoro]{Yuuki Tadokoro}
\address{Department of Mathematical Sciences,
University of Tokyo, 3-8-1 Komaba, Meguro, Tokyo 153-8914, Japan}
\email{tado\char`\@ms.u-tokyo.ac.jp}


\begin{abstract}
We prove
some value of the harmonic volume for
the Klein quartic $C$
is nonzero modulo $\displaystyle \frac{1}{2}\Z$,
using special values of the generalized hypergeometric
function ${}_3F_2$.
 This result tells us
the algebraic cycle $C-C^-$ is not
algebraically equivalent to zero
in the Jacobian variety $J(C)$. 
\end{abstract}

\maketitle
\section{Introduction}
Let $X$ be a compact Riemann surface
of genus $g\geq 2$
and $J(X)$ its Jacobian variety.
By the Abel-Jacobi map $X\to J(X)$,
$X$ is embedded in $J(X)$.
The algebraic $1$-cycle $X-X^{-}$ in $J(X)$
is homologous to zero.
Here we denote by $X^{-}$ the image of $X$
under the multiplication map by $-1$.
If $X$ is hyperelliptic,
$X=X^{-}$ in $J(X)$.
For the rest of this paper, suppose $g\geq 3$.
B. Harris \cite{H-3} studied the problem
whether the cycle $X-X^-$
in $J(X)$
is algebraically equivalent to zero or not.
The harmonic volume $I$ for $X$ was introduced by 
Harris \cite{H-1}, using Chen's iterated integrals \cite{C}.
Let
$H$ denote the first integral homology group $H_1(X; \Z)$
of $X$.
The harmonic volume $I$
is defined to be
a homomorphism $(H^{\otimes 3})^\prime\to \R/{\Z}$.
Here $(H^{\otimes 3})^\prime$
is a certain subgroup of $H^{\otimes 3}$.
See Section \ref{The harmonic volumes} for the definition.
Let $\omega$ be a third tensor product of
holomorphic $1$-forms on $X$.
Suppose that
$\omega+\overline{\omega}$
and $(\omega-\overline{\omega})/{\sqrt{-1}}$
belong to $(H^{\otimes 3})^\prime$.
If the cycle $X-X^-$ is algebraically equivalent to zero,
then
twice the values at
both $\omega+\overline{\omega}$
and $(\omega-\overline{\omega})/{\sqrt{-1}}$
of the harmonic volume
are zero modulo $\Z$.
Harris
proved twice the value at
$\omega+\overline{\omega}$
of the harmonic volume for
the Fermat quartic $F(4)$
are nonzero modulo $\Z$.
This implies
$F(4)-F(4)^-$ is not algebraically equivalent to zero
in $J(F(4))$
(\cite{H-3}, \cite{H-4}).
Ceresa \cite{Ce} showed that 
$X-X^{-}$ is not algebraically equivalent to zero
for a generic $X$.
We know few explicit nontrivial examples
except for $F(4)$.
Let $C$ denote the Klein quartic.
See Section \ref{homology basis}
for the definition.
The aim of this paper is to show
\begin{kleinthm}
The algebraic cycle $C-C^-$ is not algebraically
equivalent to zero in the Jacobian variety $J(C)$.
\end{kleinthm}
Since
Harris used the special
feature of $F(4)$
that its normalized period matrix
has entries in $\Z[\sqrt{-1}]$,
it is not difficult to find some $\omega$
so that
$\omega+\overline{\omega}$
and $(\omega-\overline{\omega})/{\sqrt{-1}}$
belong to $(H^{\otimes 3})^\prime$
for $F(4)$.
But, in general,
it is not easy to find such an $\omega$.
For the Klein quartic $C$,
we prove
$(D+\overline{D})/{7}$ and $(D-\overline{D})/{\sqrt{-7}}$
belong to $(H^{\otimes 3})^\prime$
(Proposition \ref{program}).
See Section
\ref{Some values of the harmonic volume for the Klein quartic}
for the definitions of them.
In Theorem \ref{a harmonic volume of Klein quartic}
we compute the value at
$(D-\overline{D})/{\sqrt{-7}}\in (H^{\otimes 3})^\prime$
of the harmonic volume for $C$
$$
I((D-\overline{D})/{\sqrt{-7}})=
\frac{28}{\sqrt{-7}}\biggl(
\frac{\zeta_7^{2}-\zeta_7^{6}}{\zeta_7+1}x_{1,2}
+\frac{\zeta_7^{4}-\zeta_7^{5}}{\zeta_7^{2}+1}x_{2,3}
+\frac{\zeta_7-\zeta_7^{3}}{\zeta_7^{4}+1}x_{3,1}\biggr)
\ \mathrm{mod}\ \Z.
$$
Here, $\zeta_7=\exp(2\pi\sqrt{-1}/{7})$
and $x_{i,j}$'s are real constants obtained from
some special values of the generalized hypergeometric
function ${}_3F_2$ (Lemma \ref{iterated integral computation}).
By numerical computation using MATHEMATICA,
we obtain Theorem \ref{Klein}.
We give a calculation program in Appendix.

\noindent
{\bf Acknowledgments.}
The author is grateful to Nariya Kawazumi
for valuable advice
and reading the manuscript.
Masahiko Yoshinaga and Shuji Yamamoto
suggest useful ideas for the proof of Proposition \ref{program}
to him.
He would like to thank Masaaki Suzuki
for his helpful comments for MATHEMATICA programs.
This work is partially supported by 21st
Century COE program  (University of Tokyo) by
the Ministry of Education, Culture, Sports, Science and Technology.

%\tableofcontents

\section{The harmonic volume}\label{The harmonic volumes}
\label{Preliminaries1}
We recall the harmonic volume for
a compact Riemann surface $X$ of genus $g \geq 3$ \cite{H-1}. 
We identify
the first integral homology group $H_1(X;\mathbb{Z})$
of $X$ with the first integral cohomology group
by Poincar\'e duality,
and denote it by $H$.
Moreover we identify $H$
with the space of all the
real harmonic $1$-forms on $X$
with integral periods.
Let $K$ be the kernel of the intersection pairing
$(\ , \ )\colon H\otimes_{\Z} H \to \Z$.
For the rest of this paper, we write $\otimes=\otimes_{\Z}$,
unless otherwise stated.
The Hodge star operator $\ast$
on the space of all the $1$-forms $A^1(X)$
is locally given by 
$\ast (f_1(z)dz + f_2(z)d\bar{z})
=-\sqrt{-1}f_1(z)dz + \sqrt{-1}f_2(z)d\bar{z}$
 in a local coordinate $z$
 and depends only on the complex structure
and not on the choice of Hermitian metric.
For any $\sum_{i=1}^{n}a_i\otimes b_i\in K$,
there exists a unique $\eta\in A^1(X)$
such that $d\eta=\sum_{i=1}^{n}a_i\wedge b_i$
and $\displaystyle \int_{X}\eta\wedge \ast\alpha=0$
for any closed $1$-form $\alpha\in A^{1}(X)$.
Here $a_i$ and $b_i$ are regarded as real harmonic
$1$-forms on $X$.
Choose a point $x_0\in X$.
\begin{defn}\label{pointed harmonic volume}
\mbox{(The pointed harmonic volume \cite{P})}\\
For $\sum_{i=1}^{n}a_i\otimes b_i\in K$
and $c\in H$,
the pointed harmonic volume 
defined to be
$$I_{x_0}{\Biggl(}{\biggl(}\sum_{i=1}^{n}a_{i}\otimes b_{i}{\biggr)}\otimes c{\Biggr)}=\sum_{i=1}^{n}\int_{\gamma}a_{i}b_{i}-\int_{\gamma}\eta
 \quad \mathrm{mod} \ \mathbb{Z}.$$
Here $\eta\in A^1(X)$ is associated to 
$\sum_{i=1}^{n}a_i\otimes b_i$
in the way stated above and
$\gamma$ is a loop in $X$ with the base point $x_0$ whose homology class is
equal to $c$.
The integral $\displaystyle \int_{\gamma}a_ib_i$
is Chen's iterated integral \cite{C}, that is,
$\displaystyle \int_{\gamma}a_ib_i
=\int_{0\leq t_1\leq t_2\leq 1}f_i(t_1)g_i(t_2)dt_1dt_2$
for $\gamma^{\ast}a_i=f_i(t)dt$ and $\gamma^{\ast}b_i=g_i(t)dt$.
Here $t$ is the coordinate in the interval $[0,1]$.
\end{defn}

The harmonic volume is given as a restriction of
the pointed harmonic volume $I_{x_0}$.
We denote by $(H^{\otimes 3})^{\prime}$ the kernel of
a natural homomorphism 
$p\colon H^{\otimes 3} \to H^{\oplus 3}$
defined by
$p(a\otimes b\otimes c)=((a, b)c, (b, c)a, (c, a)b)$.
The {\it harmonic volume} $I$ for $X$
is a linear form on $(H^{\otimes 3})^{\prime}$
with values in $\R/{\Z}$ defined by
the restriction of 
$I_{x_0}$ to $(H^{\otimes 3})^{\prime}$, i.e.,
$I=I_{x_0}|_{(H^{\otimes 3})^{\prime}}$.
Harris \cite{H-1}
proved that the harmonic volume $I$
is independent of the choice of the base point $x_0$.
We have 
$I(\sum_i h_{\sigma(1), i}\otimes h_{\sigma(2), i}\otimes h_{\sigma(3), i})
=\mathrm{sgn}(\sigma)I(\sum_i h_{1, i}\otimes h_{2, i}\otimes h_{3, i})\ 
\mathrm{mod}\ \Z$, 
where $\sum_i h_{1, i}\otimes h_{2, i}\otimes h_{3, i}
 \in (H^{\otimes 3})^{\prime}$
 and $\sigma$ is an element of the third symmetric group $S_3$.
See Harris \cite{H-1} and Pulte \cite{P} for details.

In general, it is difficult to compute 
the correction term $\eta$
in Definition \ref{pointed harmonic volume}.
If $X$ is a hyperelliptic curve, 
we have an explicit formula for the $1$-form $\eta$
given by Harris \cite{H-1}.
This allows us to calculate the harmonic volumes for
all the hyperelliptic curves (Tadokoro \cite{T}).
In this paper, we deal with the case $\eta$ vanishes.

\section{The algebraic cycle $X-X^{-}$ and an intermediate Jacobian}
\label{Preliminaries2}
We review a relation between
the algebraic cycle $X-X^{-}$ and the harmonic volume $I$.

Let $j_2\colon H^{\otimes 3}\to \wedge^3 H$ be
a natural homomorphism
$j_2(a\otimes b\otimes c)= 
 a\wedge b\wedge c$,
where $\wedge^3 H$ denotes the third exterior power
of $H$.
We have the homomorphism of short exact sequences
$$
\xymatrix{
0\ar[r] &(H^{\otimes 3})^{\prime} \ar[r] \ar[d]_{j_1}& 
H^{\otimes 3} \ar[d]_{j_2}\ar[r]^{p}& H^{\oplus 3} \ar[d]_{j_3}\ar[r] & 0 \\
0\ar[r] &(\wedge^3 H)^{\prime} \ar[r] & \wedge^3 H \ar[r]^{\bar{p}}& 
H \ar[r] & 0 ,
}
$$
where $j_3(a,b,c)=a+b+c$,
$\bar{p}(a\wedge b\wedge c)=(a, b)c+(b, c)a+(c, a)b$
and $j_1$ is the restriction homomorphism of $j_2$
to $(H^{\otimes 3})^{\prime}$.
Let $\mathcal{A}^{k}_{0}(J)$
be the space of algebraic
 $k$-cycles homologous to zero on the Jacobian variety
$J=J(X)$, modulo rational
equivalence.
The Abel-Jacobi map of Griffiths $\Phi_{\R}\colon 
\mathcal{A}^{k}_{0}(J)\to 
\Hom_{\Z}(H^{2k+1}(J;\Z),\R/{\Z})$
is defined by 
$$\partial W \mapsto \biggl\{ \omega\mapsto 
\int_{W}\omega\biggr\},$$
where $\omega$ is a harmonic $(2k+1)$-form 
on $J$ with integral periods
(Section 4 in \cite{P}).
Here, the module $\Hom_{\Z}(H^{2k+1}(J;\Z),\R/{\Z})$
can be identified with
an intermediate Jacobian of $H_{2k+1}(J;\Z)$ \cite{P}.
From now on, we consider the case $k=1$.
Let $\nu$ denote the Abel-Jacobi image $\Phi_{\R}(X-X^{-})$.
Harris (Proposition 2.1 in \cite{H-4}, \cite{H-1}) proved that 
$(\wedge^3 H)^{\prime}$ can be identified with
the primitive subgroup of $H^3(J;\Z)$
in the sence of Lefchetz,
denoted by $H^3_{\mathrm{prim}}(J;\Z)$.
Using this identification and the natural projection
$\Hom_{\Z}(H^{3}(J;\Z),\R/{\Z})\to 
\Hom_{\Z}(H^3_{\mathrm{prim}}(J;\Z),\R/{\Z})$,
we consider $\nu$
as an element of
$\Hom_{\Z}((\wedge^3 H)^{\prime},\R/{\Z})$
(Section 4 and 6 in \cite{P}).
\begin{thm}\label{Abel-Jacobi}$( \mathrm{Harris}$ \cite{H-1}, \cite{H-4}$)$.
The Abel-Jacobi image
$\nu$ satisfies the commutative diagram
$$
\xymatrix{
(H^{\otimes 3})^{\prime}\ar[r]^{2I} \ar[d]_{j_1}
& \R/{\Z}\\
(\wedge^3 H)^{\prime} \ar[ur]_{\nu}.&
}
$$
\end{thm}
We say the
algebraic cycle $X-X^{-}$ is
{\it algebraically equivalent to zero in} $J$
if there exists a topological $3$-chain
$W$ such that $\partial W=X-X^-$
and $W$ lies on $S$,
where $S$
is an algebraic (or complex analytic)
subset of $J$
of complex dimension $2$ (Harris \cite{H-4}).
The chain $W$ is unique up to
$3$-cycles.
We denote by $H^{1,0}$
the space of all the
holomorphc $1$-forms on $X$.
From \cite{H-3}, 2.6 in \cite{H-4} and
533-534 in \cite{We}, we have
\begin{prop}\label{cycle and intermediate Jacobian} 
Let $\omega \in \big(H^{1,0}\big)^{\otimes_{\C} 3}$
satisfying that 
$\omega+\overline{\omega}$ and $(\omega-\overline{\omega})/{\sqrt{-1}}
\in (H^{\otimes 3})^{\prime}$.
If $X-X^-$ is algebraically
equivalent to zero in $J$,
then twice the values at both
$\omega+\overline{\omega}$ 
and $(\omega-\overline{\omega})/{\sqrt{-1}}$
of the harmonic volume are zero
modulo $\Z$.
\end{prop}
\begin{proof}
Since $X-X^-$ is algebraically
equivalent to zero in $J$,
there exist a $3$-chain $W$ and
an algebraic subset $S$
satisfying the above conditions.
Let $H_{\C}$ denote 
$H\otimes \C$.
Theorem
\ref{Abel-Jacobi}
gives
$$2I(\omega+\overline{\omega})=
 \int_{W}j_1(\omega+\overline{\omega})\ \textnormal{and}\ 
2I((\omega-\overline{\omega})/{\sqrt{-1}}))=
 \int_{W}j_1(\omega-\overline{\omega})/{\sqrt{-1}}.$$
It is clear that 
$j_1(\omega)$ and $j_1(\overline{\omega})$
are $(3,0)$ and $(0,3)$-form 
in $H^3(J;\C)=\wedge^3 H_{\C}$ respectively.
Since $\dim_{\C}S=2$,
the restriction of
them to $S$ are clearly zero.
\end{proof}
If twice the value at
$\omega+\overline{\omega}$ 
or $(\omega-\overline{\omega})/{\sqrt{-1}}$
of the harmonic volume
is nonzero modulo $\Z$,
then $X-X^-$ is not algebraically equivalent to zero in $J$.
See Hain \cite{Hain}, Pirola \cite{Pirola}
and their references for
the algebraic cycle $X-X^{-}$ in $J$.

\section{Some values of the harmonic volume for the Klein quartic}
We compute some values of the harmonic volume for the Klein quartic
to prove the main theorem (Theorem \ref{Klein}).
%Klein curve 'Ì'è‹`
\subsection{A $1$-dimensional homology basis of the Klein quartic}
\label{homology basis}
We denote by $C$ the Klein quartic
which is, by definition,
the plane curve
$C\colon=\{(X:Y:Z)\in\C P^2; X^3Y+Y^3Z+Z^3X=0\}$.
It is a compact Riemann surface of genus 3.
It is known that the holomorphic automorphism group of $C$,
$\Aut(C)$, is isomorphic to 
$\mathrm{PSL}_2(\mathbb{F}_7)$.
See \cite{S} for the details
of the Klein quartic. 
Let $x$ and $y$ denote $X^3Y^{-2}Z^{-1}+1$ and 
$-XY^{-1}$ respectively. 
The equation $X^3Y+Y^3Z+Z^3X=0$
induces $y^7=x(1-x)^2$.
The holomorphic map $\pi: C\to\C P^1$ is defined by
 $\pi(x, y)=x$,
which is a $7$-sheeted covering $C \to \C P^1$, 
branched over $3$ branch points $\{0,1,\infty\}$.
Let $\zeta_7$ denote $\exp(2\pi\sqrt{-1}/{7})$.
For $t\in [0,1]$, we define a loop
 $e_0: [0,1]\to C$ by
$e_0(t)=(t,y_0(t))$,
where $y_0(t)$ is a real analytic
function $\sqrt[7]{t(1-t)^2}$.
Let $\sigma\colon C\to C$ be a holomorphic automorphism
$\sigma(x,y)=(x,\zeta_7 y)$.
For $k=0,1,\ldots, 6$, we define loops in $C$
by $c_k=\sigma_{\ast}^{k}(e_0)\cdot e_0^{-1}$.
We denote 
$\ell_k=\sigma_\ast^{k-1}(e_0)\cdot 
\sigma_\ast^{k}(e_0)^{-1}, k=0,1,\ldots,7$.
The loop $\ell_0$ can be identified 
with $\ell_7$.
By abuse of notation, the homology classes of $c_k$ and $\ell_k$
are denoted by $c_k$ and $\ell_k\in H_1(C;\Z)$
respectively.
Let $(\ ,\ ): H_1(C;\Z)\otimes H_1(C;\Z)\to \Z$
 be the intersection pairing, i.e., a non-degenerate
bilinear form on $H_1(C;\Z)$.
Tretkoff and Tretkoff \cite{T-T}
proved
$$
(c_1, c_k)=
\left\{
\begin{array}{crl}
0&\mathrm{if} & k=1,2,4,6, \\
1&\mathrm{if} & k=3,5,
\end{array}
\right.
$$
using the Hurwitz system
of the branched covering $\pi$.
By the definition of $\ell_k$,
we have
$$
(\ell_1, \ell_k)=(c_1,c_k)
-(c_1,c_{k-1})=
\left\{
\begin{array}{crl}
0&\mathrm{if} & k=1,2, \\
1&\mathrm{if} & k=3,5,\\
-1&\mathrm{if} & k=4,6.
\end{array}
\right.
$$
Moreover,
we obtain that 
$\sigma_{\ast}(\ell_k)=\ell_{k+1}$
and $(\ell_i,\ell_j)=
(\sigma_{\ast}(\ell_{i}), \sigma_{\ast}(\ell_{j}))=
(\ell_{i+1},\ell_{j+1})$.
The intersection matrix $K^{\prime}$
of $\ell_k, k=1,2,\ldots,6$
is given by
$$
\left(
\begin{array}{cccccc}
0 &0 &1 &-1 &1 &-1 \\
0 &0 &0 &1 & -1&1 \\
-1 &0 &0 &0 &1 &-1 \\
1 &-1 &0 &0 &0 &1 \\
-1 &1 &-1 &0 &0 &0 \\
1 &-1 &1 &-1 &0 &0 \\
\end{array}
\right),
$$
i.e.,
its $(i,j)$-th entry is $(\ell_i,\ell_j)$.
It is easy to prove $\det K^{\prime}=1$ and 
$\{\ell_k\}_{k=1,2,\ldots,6}\subset H_1(C;\Z)$
is a basis of $H_1(C;\Z)$.

\subsection{Poincar\'e dual of the Klein quartic}
Let $\omega^{\prime}_1, \omega^{\prime}_2$
and $\omega^{\prime}_3$ be
 holomorphic $1$-forms on $C$, $(1-x)dx/{y^6},
(1-x)dx/{y^5}$ and $dx/{y^3}$ respectively.
It is known that
$\{\omega^{\prime}_i\}_{i=1,2,3}$ is a basis of
the space of all the holomorphic $1$-forms on $C$.
The beta function $B(u,v)$ is defined by $\displaystyle\int_0^1t^{u-1}
(1-t)^{v-1}dt$ for $u,v>0$. 
We denote $(h_1,h_2,h_3,h_4)=(1/{7},2/{7},4/{7},1/{7})$
and $\xi_i=\zeta_7^{7h_i}$.
From the equations
$\sigma^{\ast}\omega^{\prime}_i=\xi_i\omega_i$
and
$\displaystyle\int_{e_0}\omega^{\prime}_i=B(h_i,h_{i+1})$,
we have
\begin{lem}\label{period of Klein quartic}
$$
\int_{\ell_k}\omega^{\prime}_i=
(\xi_i^{k-1}-\xi_i^k)B(h_i,h_{i+1}).
$$
\end{lem}
\begin{rem}
These integrals depend only on the cohomology class 
of $\omega_j^{\prime}$
 and the homology class of $\ell_k$.
\end{rem}

We set
$B^{\prime}_i=B(h_i,h_{i+1})$
and
$\omega_i=
\omega^{\prime}_i/{B^{\prime}_i}$, $i=1,2,3$.
We write $L_k:=\sum_{i=1}^{7}\zeta_7^{ik}\ell_k\in H_1(C;\C)$
and
denote the Poincar\'e dual
by $\mathrm{P.D.}\colon H^1(C;\C)\to H_1(C;\C)$.
\begin{prop}\label{dual}
We denote $\lambda_i=-1/(\xi_i^3(\xi_i^2+1))\in \C$.
Then, we have
$$
\mathrm{P.D.}(\omega_{i})=
\lambda_iL_{7h_i}.
$$
\end{prop}
\begin{proof}
Since $\sigma_{\ast}(\ell_k)=\ell_{k+1}$,
we obtain 
$\sigma_{\ast}L_k=\zeta_7^{-k}L_k$.
The eigenvalues and eigenvectors of
the action of $\sigma$ on the 
$\C$-vector space $H_1(C;\C)$
are $\zeta_7^{-k}$ and $L_k$ for $k=1,2,\ldots,6$.
We have
$$
\sigma_{\ast}(\mathrm{P.D.}(\omega_{i}))
=\mathrm{P.D.}((\sigma^{-1})^{\ast}\omega_{i})
=\xi_i^{-1}\mathrm{P.D.}(\omega_{i})
=\zeta_7^{-7h_i}\mathrm{P.D.}(\omega_{i}).
$$
There exists a constant $\lambda_i\in \C$ such that
$\mathrm{P.D.}(\omega_{i})=\lambda_iL_{7h_i}$.
The result follows from Lemma \ref{period of Klein quartic}
and the equation
$$
\int_{\ell_1}\omega_i=(\mathrm{P.D.}(\omega_{i}),\ell_1)
=(\lambda_iL_{7h_i},\ell_1)=\lambda_i(L_{7h_i},\ell_1)
=-\lambda_i(1-\xi_i)(\xi_i^3(\xi_i^2+1)).
$$
\end{proof}
\begin{rem}\label{rmk}
We have $\mathrm{P.D.}(\overline{\omega}_{i})=
\overline{\lambda}_i\overline{L}_{7h_i}$.
It immediately follows $\lambda_1\lambda_2\lambda_3=-1$.
\end{rem}

\subsection{Some values of the harmonic volume for the Klein quartic}
\label{Some values of the harmonic volume for the Klein quartic}
For $t\in[0,1]$,
let $f_i$ be
a real $1$-form on $[0,1]$
defined by
$e_0^{\ast}\omega^{\prime}_i
=t^{h_{i}-1}(1-t)^{h_{i+1}-1}dt$, $i=1,2,3$.
Let $x_{i,j}$ denote an iterated integral 
$\displaystyle \int_{e_0}\omega_i\omega_j=
\int_{\gamma}f_if_j \Big/{(B^{\prime}_iB^{\prime}_j)}$.
Here, $\gamma$ is the path $[0,1]\ni t\mapsto
t\in[0,1]$.
We compute the iterated integrals of $\omega_1,\omega_2$ and
$\omega_3$ along the loop $\ell_k$.
\begin{lem}\label{iterated integral of Klein quartic}
We consider $\ell_k$ as loops with the base point
$(x,y)=(0,0)\in C$. We have 
$$
\int_{\ell_k}\omega_i\omega_j=
(\xi_i\xi_j)^{k-1}(1-\xi_i\xi_j)x_{i,j}
+(\xi_i\xi_j)^{k-1}(\xi_i\xi_j-\xi_j).
$$
\end{lem}
\begin{rem}
Since $\omega_i$ is closed and $\omega_i\wedge \omega_j=0$,
these iterated integrals are invariant under homotopy with fixed
endpoints.
\end{rem}
\begin{proof}
Using 
the shuffle product formula (Chen \cite{C}, 1.6)
and the equations
$$0=\int_{e_0\cdot e_0^{-1}}\omega_i\omega_j
=\int_{e_0}\omega_i\omega_j
+\int_{e_0^{-1}}\omega_i\omega_j
+\int_{e_0}\omega_i\int_{e_0^{-1}}\omega_j
\quad \textnormal{and}\quad
\int_{e_0}\omega_i=\int_{e_0}\omega^{\prime}_i
\biggl/{B_i^{\prime}}=1,$$
we have
\begin{align*}
\int_{\ell_k}\omega_i\omega_j=&
\int_{\sigma_\ast^{k-1}(e_0)\cdot 
\sigma_\ast^{k}(e_0)^{-1}
}\omega_i\omega_j\\
=&\int_{\sigma_\ast^{k-1}(e_0)}\omega_i\omega_j
+\int_{\sigma_\ast^{k}(e_0)^{-1}}\omega_i\omega_j
+\int_{\sigma_\ast^{k-1}(e_0)}\omega_i
\int_{\sigma_\ast^{k}(e_0)^{-1}}\omega_j\\
=&(\xi_i\xi_j)^{k-1}\int_{e_0}\omega_i\omega_j
+(\xi_i\xi_j)^{k}\int_{e_0^{-1}}\omega_i\omega_j
-\xi_i^{k-1}\xi_j^{k}\int_{e_0}\omega_i\int_{e_0}\omega_j\\
=&(\xi_i\xi_j)^{k-1}\int_{e_0}\omega_i\omega_j
+(\xi_i\xi_j)^{k}
\biggl\{
-\int_{e_0}\omega_i\omega_j
+\int_{e_0}\omega_i\int_{e_0}\omega_j
\biggr\}
-\xi_i^{k-1}\xi_j^{k}\\
=&(\xi_i\xi_j)^{k-1}(1-\xi_i\xi_j)\int_{e_0}\omega_i\omega_j
+(\xi_i\xi_j)^{k-1}(\xi_i\xi_j-\xi_j).
\end{align*}
\end{proof}

The subset $\mathcal{H}$
of $(H^{\otimes 3})^{\prime}\otimes\R$
is defined by
$\{\omega+\overline{\omega},
(\omega-\overline{\omega})/{\sqrt{-1}};
\omega\in H^{1,0}\otimes_{\C} H^{1,0}\otimes_{\C} H^{1,0}\}$.
We will find some elements of
$\mathcal{H}\cap (H^{\otimes 3})^{\prime}$.
Let $D$ and $\overline{D}$ denote
$\sum_{\mu\in S_3}\mathrm{sgn}(\mu)\omega_{\mu(1)}
\otimes_{\C} \omega_{\mu(2)}\otimes_{\C} \omega_{\mu(3)}$
and
$\sum_{\mu\in S_3}\mathrm{sgn}(\mu)
\overline{\omega}_{\mu(1)}\otimes_{\C}
\overline{\omega}_{\mu(2)}\otimes_{\C}\overline{\omega}_{\mu(3)}
\in (H_{\C})^{\otimes_{\C} 3}$
respectively.
Using 
Proposition \ref{dual}
and 
Remark \ref{rmk},
$D$ and $\overline{D}$ are identified with
$-\sum_{\mu\in S_3}\mathrm{sgn}(\mu)
L_{7h_{\mu(1)}}\otimes_{\C} L_{7h_{\mu(2)}}\otimes_{\C} L_{7h_{\mu(3)}}$
and
$-\sum_{\mu\in S_3}\mathrm{sgn}(\mu)
\overline{L}_{7h_{\mu(1)}}\otimes_{\C} \overline{L}_{7h_{\mu(2)}}
\otimes_{\C} \overline{L}_{7h_{\mu(3)}}$
respectively.
The coefficients of 
$\ell_p\otimes_{\C} \ell_q\otimes_{\C} \ell_r$
of $D$ and $\overline{D}$
are 
$$
\alpha_{p,q,r}=-
\left|
\begin{array}{ccc}
\zeta_7^{p} & \zeta_7^{2p} & \zeta_7^{4p}\\
\zeta_7^{q} & \zeta_7^{2q} & \zeta_7^{4q}\\
\zeta_7^{r} & \zeta_7^{2r} & \zeta_7^{4r}
\end{array}
\right|
\ \textnormal{and}\quad 
\overline{\alpha}_{p,q,r}=-
\left|
\begin{array}{ccc}
\zeta_7^{6p} & \zeta_7^{5p} & \zeta_7^{3p}\\
\zeta_7^{6q} & \zeta_7^{5q} & \zeta_7^{3q}\\
\zeta_7^{6r} & \zeta_7^{5r} & \zeta_7^{3r}
\end{array}
\right|
$$
respectively.
It is trivial that 
$D+\overline{D}$ and $(D-\overline{D})/{\sqrt{-1}}\in \mathcal{H}$.
Furthermore, we have
\begin{prop}\label{program}
$(D+\overline{D})/{7}$ and $(D-\overline{D})/{\sqrt{-7}}
\in (H^{\otimes 3})^\prime$.
\end{prop}
\begin{proof}
It suffices to prove that
$\alpha_{p,q,r}$ belongs to the principal ideal 
$(\sqrt{-7})\Z[(1+\sqrt{-7})/{2}]
\subset\Z[(1+\sqrt{-7})/{2}]$.
It is well known that
$\Gal(\Q(\zeta_7)/{\Q})\cong\{\sigma_i\}_{i=1,2,\ldots,6}
\cong\Z/{6\Z}$,
where $\sigma_i(\zeta_7)=\zeta_7^{i}$.
Since $[\Q(\sqrt{-7})\colon \Q ]=2$,
we obtain 
$\Gal(\Q(\zeta_7)/{\Q(\sqrt{-7})})$,
the subgroup of $\Gal(\Q(\zeta_7)/{\Q})$,
is generated by $\sigma_2$.
It is clear that
$\alpha_{p,q,r}$ is invariant under the action of $\sigma_2$.
So, we have $\alpha_{p,q,r}\in \Q(\sqrt{-7})$.
On the other hand,
it immediately follows 
$\alpha_{p,q,r}$ belongs to 
the principal ideal $(\zeta_7-1)\Z[\zeta_7]\subset \Z[\zeta_7]$.
Therefore, we have 
$$\alpha_{p,q,r}\in \Q(\sqrt{-7}) \cap 
(\zeta_7-1)\Z[\zeta_7]
=(\sqrt{-7})\Z[(1+\sqrt{-7})/{2}]\subset 
\Z[(1+\sqrt{-7})/{2}].$$
We have
$\alpha_{p,q,r}+\overline{\alpha}_{p,q,r}
\in 7\Z$ and 
$\alpha_{p,q,r}-\overline{\alpha}_{p,q,r}\in \sqrt{-7}\Z$.
We complete the proof.
\end{proof}
\begin{rem}
Using the character of $\mathrm{Aut}(C)
=PSL_2(\mathbb{F}_7)$, we have 
$H^0(\mathrm{Aut}(C);(H_{\C})^{\otimes_{\C} 3})
=\C^{2}$.
This induces
$H^0(\mathrm{Aut}(C);H^{\otimes 3})
=\Z^{2}$.
We can also prove that
$\{(D+\overline{D})/{7},(D-\overline{D})/{\sqrt{-7}}\}$
is a generator of
 $H^0(\mathrm{Aut}(C);(H^{\otimes 3})^\prime)$.
\end{rem}
\begin{thm}\label{a harmonic volume of Klein quartic}
The values at
$(D+\overline{D})/7$ and $(D-\overline{D})/{\sqrt{-7}}
\in (H^{\otimes 3})^\prime$
for the harmonic volume of the Klein quartic $C$
are given by
$$
0
\ \textit{and}\ 
\frac{28}{\sqrt{-7}}\biggl(
\frac{\zeta_7^{2}-\zeta_7^{6}}{\zeta_7+1}x_{1,2}
+\frac{\zeta_7^{4}-\zeta_7^{5}}{\zeta_7^{2}+1}x_{2,3}
+\frac{\zeta_7-\zeta_7^{3}}{\zeta_7^{4}+1}x_{3,1}\biggr)
\ \mathrm{mod}\ \Z$$
respectively.
\end{thm}
\begin{proof}
All iterated integral parts of
$I((D+\overline{D})/7)$ and $I((D-\overline{D})/{\sqrt{-7}})$
are linear combinations of
$\displaystyle\int_{\ell_k}\omega_i\omega_j$ and
$\displaystyle\int_{\ell_k}\overline{\omega}_i\overline{\omega}_j
=\overline{\int_{\ell_k}\omega_i\omega_j}$.
Furthermore, $\omega_i\wedge \omega_j=\overline{\omega}_i\wedge
\overline{\omega}_j=0$.
So we need no correction terms $\eta$
in Definition \ref{pointed harmonic volume}.
Therefore, it suffices to
calculate only the iterated integral parts.

By definition,
there exist complex constants $\theta_{i,j,k}$
so that
$I((D+\overline{D})/{7})$
is of the form
$$\sum_{k=1}^{7}\sum_{(i,j)\in U}
\theta_{i,j,k}\int_{\ell_k}(\omega_i\omega_j-
\omega_j\omega_i)
+\sum_{k=1}^{7}\sum_{(i,j)\in U}
\overline{\theta}_{i,j,k}
\overline{\int_{\ell_k}(\omega_i\omega_j-
\omega_j\omega_i)},$$
where $U$ is a set $\{(1,2),(2,3),(3,1)\}$.
Using $\mathrm{P.D.}(\omega_{i})=\lambda_iL_{7h_i}
=\lambda_i\sum_{k=1}^{7}\xi_i^k\ell_k$,
it can be written as
$(I_{1,2,3}+\overline{I}_{1,2,3})/{7}
\ \mathrm{mod}\ \Z$.
Here, we denote
$$
I_{1,2,3}=\lambda_3\sum_{k=1}^{7}\xi_3^{k}\int_{\ell_k}
(\omega_1\omega_2-\omega_2\omega_1)
+\lambda_1\sum_{k=1}^{7}\xi_1^k\int_{\ell_k}
(\omega_2\omega_3-\omega_3\omega_2)
+\lambda_2\sum_{k=1}^{7}\xi_2^k\int_{\ell_k}
(\omega_3\omega_1-\omega_1\omega_3).
$$
Similarly, we obtain
$$I((D-\overline{D})/{\sqrt{-7}})=
(I_{1,2,3}-\overline{I}_{1,2,3})/{\sqrt{-7}}
\ \mathrm{mod}\ \Z.$$
In order to complete the proof,
we need two lemmas.
\begin{lem}\label{lemma for proposition}
We have
$$
\int_{\ell_k}(\omega_i\omega_j-\omega_j\omega_i)
=2(\xi_i\xi_j)^{k-1}(1-\xi_i\xi_j)x_{i,j}
+(\xi_i\xi_j)^{k-1}(\xi_i-1)(\xi_j+1).
$$
\end{lem}
\begin{proof}
We use 
Lemma \ref{period of Klein quartic},
Lemma \ref{iterated integral of Klein quartic}
and the equation
$$
\int_{\ell_k}\omega_j\omega_i
=-\int_{\ell_k}\omega_i\omega_j
+\int_{\ell_k}\omega_i\int_{\ell_k}\omega_j.
$$
\end{proof}
\begin{lem}\label{iterated integral sum}
We have
$$I_{1,2,3}=14\biggl(
\frac{\zeta_7^{2}-\zeta_7^{6}}{\zeta_7+1}x_{1,2}
+\frac{\zeta_7^{4}-\zeta_7^{5}}{\zeta_7^{2}+1}x_{2,3}
+\frac{\zeta_7-\zeta_7^{3}}{\zeta_7^{4}+1}x_{3,1}
-\frac{3}{2}\sqrt{-7}
\biggr).$$
\end{lem}
\begin{proof}
Using Lemma \ref{lemma for proposition}
and $\xi_1\xi_2\xi_3=1$,
we calculate the coefficient of $x_{1,2}$ of $I_{1,2,3}$
as follows:
\begin{align*}
\lambda_3\sum_{k=1}^{7}\xi_3^k\cdot 2(\xi_1\xi_2)^{k-1}(1-\xi_1\xi_2)
=&\frac{-2}{\xi_3^3(\xi_3^2+1)}
\sum_{k=1}^{7}(\xi_1\xi_2\xi_3)^{k-1}(\xi_3-1)\\
=&\frac{-2}{\zeta_7^{12}(\zeta_7^{8}+1)}
\sum_{k=1}^{7}(\zeta_7^{4}-1)
=14\frac{\zeta_7^{2}-\zeta_7^{6}}{\zeta_7+1}.
\end{align*}
Similarly, we compute
the coefficients of $x_{2,3}, x_{3,1}$ and
the constant term of $I_{1,2,3}$.
For the computation of the constant term,
we need $\zeta_7+\zeta_7^2+\zeta_7^4=(-1+\sqrt{-7})/{2}$.
\end{proof}
The result follows from
Lemma \ref{iterated integral sum}.
We remark that
all the coefficients of $x_{1,2}, x_{2,3}, x_{3,1}$ and
the constant term of $I_{1,2,3}$
are pure imaginary.
\end{proof}

For the numerical calculation of 
$x_{i,j}$,
we recall the generalized hypergeometric function ${}_3F_2$.
We denote the gamma function
$\Gamma(\tau)=
\displaystyle\int_{0}^{\infty}e^{-t}t^{\tau-1}dt$
for $\tau>0$
and $(\alpha, n)=\G(\alpha +n)/{\G(\alpha)}$
for non-negative integer $n$.
For $x\in\{z\in \C; |z|<1\}$
and $\alpha_1,\alpha_2,\alpha_3,\beta_1,\beta_2>-1$,
the generalized hypergeometric function ${}_3F_2$
is defined by
$${}_3F_2{\Big(}
\left.
\begin{array}{c}
\alpha_1,\alpha_2,\alpha_3\\
\beta_1,\beta_2
\end{array}
\right.
;x{\Big)}
=\sum_{n=0}^{\infty}{{(\alpha_1,n)(\alpha_2,n)(\alpha_3,n)}
\over{(\beta_1,n)(\beta_2,n)(1,n)}}x^n.
$$
See \cite{Sl} for example.
By straightforward computation, we have
\begin{prop}\label{hypergeometric function}
Let $\Delta$ be a $1$-simplex
$\{(u,v)\in\R^2;0\leq v\leq 1,
0\leq u\leq v\}$.
If $a,b,p,q>0, b<1$, then we have\\
$\displaystyle
\int_{\Delta} u^{a-1}(1-u)^{b-1}
v^{p-1}(1-v)^{q-1}dudv
={{B(a+p,q)}\over{a}}
\lim_{
\begin{subarray}{c}
t\to 1-0\\
t\in \R
\end{subarray}
}{}_3F_2{\Big(}
\left.
\begin{array}{c}
a,1-b,a+p\\
1+a,a+p+q
\end{array}
\right.
;t{\Big)}$.
\end{prop}

From Proposition \ref{hypergeometric function},
we have
\begin{lem}\label{iterated integral computation}
$$x_{i,j}={{B(h_{i}+h_{j},h_{j+1})}
\over
{h_{i}B_i^{\prime}B_j^{\prime}}}
\lim_{
\begin{subarray}{c}
t\to 1-0\\
t\in \R
\end{subarray}
}{}_3F_2{\Big(}
\left.
\begin{array}{c}
h_{i},1-h_{i+1},h_{i}+h_{j}\\
1+h_{i},h_{i}+h_{j}+h_{j+1}
\end{array}
\right.
;t{\Big)}.$$
\end{lem}

%It is clear that
%$\displaystyle\int_{\gamma}f_if_j+\int_{\gamma}f_jf_i=
%\int_{\gamma}f_i\int_{\gamma}f_j$.
%Then, we have
%\begin{equation}\label{equation}
%\int_{\gamma}f_if_j-\int_{\gamma}f_jf_i
%=2\int_{\gamma}f_if_j-\int_{\gamma}f_i\int_{\gamma}f_j.
%\end{equation}
\begin{thm}\label{Klein}
Let $C$ be the Klein quartic. Then,
the cycle $C-C^-$ is not algebraically
equivalent to zero in $J(C)$.
\end{thm}
\begin{proof}
By
Theorem \ref{a harmonic volume of Klein quartic},
Lemma \ref{iterated integral computation},
%and the equation (\ref{equation}),
the numerical calculation (Figure \ref{klein-program}
in Appendix),
we obtain the value
$$2I((D-\overline{D})/{\sqrt{-7}})=
0.72270\pm 1\times 10^{-5}
\ \mathrm{mod}\ \Z.$$
The result follows from
Proposition \ref{cycle and intermediate Jacobian}.
\end{proof}

\section{Appendix}
In this section, we introduce the MATHEMATICA program \cite{W} in
the proof of Theorem \ref{Klein}.
\begin{figure}[h]
\begin{center}
\psbox[width=15cm]{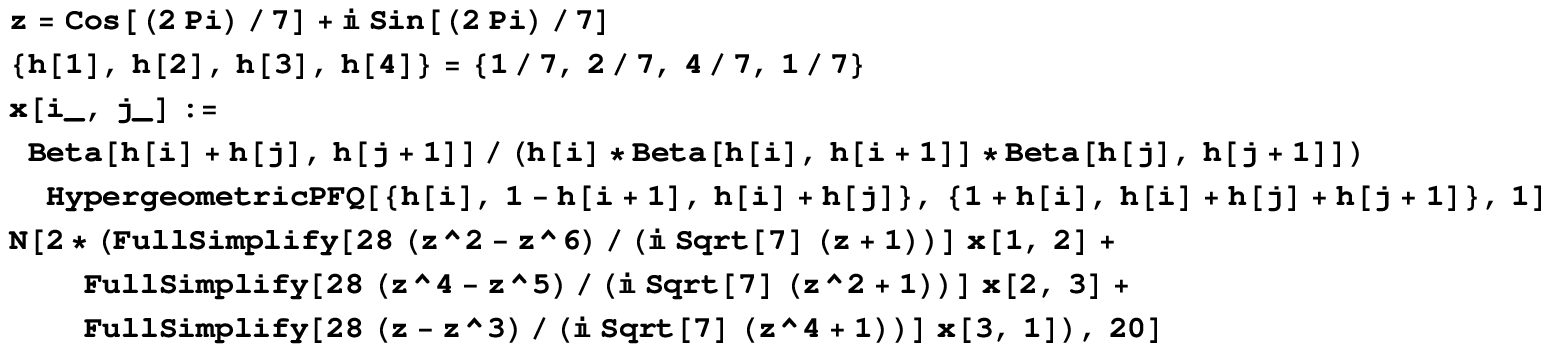}
\caption{Numerical calculation program of Theorem \ref{Klein}}
\label{klein-program}
\end{center}
\end{figure}

\end{document}